\documentclass[12pt]{article}
\usepackage[utf8]{inputenc}
\usepackage{lmodern}
\usepackage{amsmath,amssymb,amsthm}
\usepackage{booktabs}
\usepackage{geometry}
\newtheorem{theorem}{Theorem}[section]

\newtheorem{lemma}[theorem]{Lemma}

\theoremstyle{definition}
\newtheorem{problem}{Problem}[section]
\theoremstyle{remark}
\newtheorem{remark}[theorem]{Remark}

\DeclareMathOperator{\Aut}{Aut}
\DeclareMathOperator{\Inn}{Inn}
\DeclareMathOperator{\Hom}{Hom}
\DeclareMathOperator{\Syl}{Syl}
\newcommand{\Z}{\mathbb{Z}}
\newcommand{\N}{\mathbb{N}}

\title{On some problems regarding a group-theoretical generalization of the Gauss formula}
\author{Marius T\u arn\u auceanu}
\date{November 15, 2026}

\begin{document}

\maketitle

\begin{abstract}
In \cite{FT23}, we introduced the group-theoretical function $f(G) = \frac{1}{|G|}\sum_{H \le G} |\Aut(H)|$ for a finite group $G$ and posed three open problems concerning its values on noncyclic groups. In this note, we provide complete solutions to all three problems.  
\end{abstract}
\smallskip

{\small
\noindent
{\bf MSC 2020\,:} Primary 20D60, 11A25; Secondary 20D99, 11A99.

\noindent
{\bf Key words\,:} Gauss formula, Euler's totient function, automorphism group, finite group, cyclic group, $p$-group.}

\section{Introduction}

For a finite group $G$, consider the functions
\[
S(G) = \sum_{H \le G} |\Aut(H)| \quad \text{and} \quad f(G) = \frac{S(G)}{|G|}.
\]
As shown in \cite{FT23}, $f$ is multiplicative on direct products of finite groups of coprime orders, and $f(G) \ge 1$ with equality if and only if $G$ is cyclic. Denoting by $C(G)$ the poset of cyclic subgroups of $G$, we can write
\begin{equation}\label{eq:f_noncyclic}
f(G) = 1 + \frac{1}{|G|} \sum_{H \notin C(G)} |\Aut(H)|.
\end{equation}
Observe also that $|\Aut(K)| \ge 4$ for every finite noncyclic group $K$. Indeed, if $K$ is abelian and noncyclic, then $|\Aut(K)| > \varphi(|K|)$ by \cite[Corollary 2.2]{FT23} (with $|K| \ge 4$ and $|\Aut(K)|$ even), whereas if $K$ is nonabelian, then $|\Aut(K)| \ge |\Inn(K)| = |K : Z(K)| \ge 4$.

In \cite{FT23}, the following three open problems were proposed:

\begin{problem}\label{prob1}
Determine all finite groups $G$ such that $f(G) = 2$.
\end{problem}

\begin{problem}\label{prob2}
Find the minimum of $f$ on the class of finite noncyclic $p$-groups.
\end{problem}

\begin{problem}\label{prob3}
Is $\mathrm{Im}(f)$ dense in the interval $[2, \infty)$?
\end{problem}

The purpose of this note is to settle all three problems. We first solve Problem \ref{prob2} by showing that the minimum of $f$ on the class of finite noncyclic $p$-groups is $\frac{5}{2}$, attained uniquely for $p=2$ and $G \cong \Z_2 \times \Z_2$. Next, we answer Problem \ref{prob1} by proving that $f(G) = 2$ if and only if $G \cong (\Z_3 \rtimes \Z_{2^k}) \times \Z_m$, where $k \in \N^*$, the action of $\Z_{2^k}$ on $\Z_3$ is nontrivial, and $m$ is a positive integer coprime to $6$. Finally, we establish that $\mathrm{Im}(f) \cap \bigl(2, \frac{5}{2}\bigr) = \emptyset$, which answers Problem \ref{prob3} in the negative.

\section{Preliminary observations}

Before proving our main theorems, we record a brief lemma collecting elementary lower bounds that will be used repeatedly throughout our inductive arguments.

\begin{lemma}\label{lem:basic}
Let $G$ be a finite group and let $M \le G$ be a subgroup of index $r = [G : M]$.
\begin{enumerate}
    \item[\textup{(a)}] One has
    \[
    f(G) \ge 1 + \frac{1}{r}\bigl(f(M) - 1\bigr) + \frac{1}{|G|} \sum_{\substack{H \le G,\, H \not\subseteq M \\ H \notin C(G)}} |\Aut(H)|.
    \]
    \item[\textup{(b)}] If $G$ is noncyclic and $M < G$ is proper, then $f(G) > 1 + \frac{1}{r}\bigl(f(M) - 1\bigr)$. In particular, if $M$ is noncyclic, then $f(G) > 1$.
    \item[\textup{(c)}] If $G = M \times \Z_p$ for a prime $p$, then
    \[
    f(G) \ge 1 + \frac{p}{r}\bigl(f(M) - 1\bigr) = f(M).
    \]
\end{enumerate}
\end{lemma}

\begin{proof}
Part (a) follows immediately from \eqref{eq:f_noncyclic} by splitting the sum over noncyclic subgroups $H \le G$ into those contained in $M$ and those not contained in $M$, and noting that $\sum_{H \le M,\, H \notin C(M)} |\Aut(H)| = |M|\bigl(f(M) - 1\bigr)$.

Part (b) follows from (a) because $G \notin C(G)$ and $G \not\subseteq M$, with $|\Aut(G)| \ge 4$.

For part (c), note that $r = p$ and every noncyclic subgroup $H \le M$ produces a distinct noncyclic subgroup $H \times \Z_p \not\subseteq M$ of $G$ satisfying $|\Aut(H \times \Z_p)| \ge \varphi(p)|\Aut(H)| = (p-1)|\Aut(H)|$. Summing over all $H \le M$ with $H \notin C(M)$ yields
\[
\sum_{\substack{N \le G,\, N \not\subseteq M \\ N \notin C(G)}} |\Aut(N)| \ge (p-1)|M|\bigl(f(M) - 1\bigr),
\]
and substituting this estimate into (a) gives $f(G) \ge f(M)$.
\end{proof}

\section{Main results}

We begin by answering Problem \ref{prob2}.

\begin{theorem}\label{thm:pgroups}
Let $G$ be a finite noncyclic $p$-group. Then
\[
f(G) \ge \frac{5}{2},
\]
and equality holds if and only if $p = 2$ and $G \cong \Z_2 \times \Z_2$.
\end{theorem}

\begin{proof}
Let $|G| = p^n$ with $n \ge 2$. If $G$ is abelian, the conclusion follows from \cite[Theorem 1.4]{FT23}. Assume henceforth that $G$ is a nonabelian $p$-group, so $n \ge 3$. We show by induction on $n \ge 3$ that $f(G) > \frac{5}{2}$.

\medskip
\noindent\textbf{Base case ($n = 3$).} Let $G$ be a nonabelian group of order $p^3$.
\begin{itemize}
    \item If $p = 2$, then $G \cong D_8$ or $G \cong Q_8$. Direct computation gives
    \[
    f(D_8) = \frac{1 + 5 \cdot 1 + (2 + 2 \cdot 6) + 8}{8} = \frac{22}{8} = \frac{11}{4} > \frac{5}{2} \quad \text{and} \quad f(Q_8) = 4 > \frac{5}{2}\,.
    \]
    \item If $p \ge 3$, then $G$ is isomorphic either to $(\Z_p \times \Z_p) \rtimes \Z_p$, for which $|\Aut(G)| = p^3(p-1)^2(p+1)$, or to $\Z_{p^2} \rtimes \Z_p$, for which $|\Aut(G)| = p^3(p-1)^2$. In both cases, \eqref{eq:f_noncyclic} yields
    \[
    f(G) > 1 + \frac{|\Aut(G)|}{p^3} \ge 1 + (p-1)^2 \ge 5 > \frac{5}{2}\,.
    \]
\end{itemize}

\medskip
\noindent\textbf{Inductive step ($n \ge 4$).} Suppose that $n \ge 4$ and that $f(H) \ge \frac{5}{2}$ for every noncyclic $p$-group $H$ of order $p^{n-1}$. Since $G$ is a noncyclic $p$-group, it possesses at least $p + 1$ maximal subgroups $M_0, M_1, \dots, M_p$, each of order $p^{n-1}$.

If all maximal subgroups of $G$ are cyclic, then $G$ is a minimal noncyclic $p$-group; being nonabelian, \cite[Proposition 1.3]{FT23} implies $G \cong Q_8$, for which $f(Q_8) = 4 > \frac{5}{2}$.

Assume now that at least one maximal subgroup, say $M_0$, is noncyclic. By the inductive hypothesis, $S(M_0) \ge \frac{5}{2}p^{n-1}$, whence
\[
\sum_{\substack{H \le M_0 \\ H \notin C(M_0)}} |\Aut(H)| = S(M_0) - p^{n-1} \ge \frac{3}{2}p^{n-1}.
\]
Furthermore, every group $M_i$ ($1 \le i \le p$) of order $p^{n-1}$ satisfies $|\Aut(M_i)| \ge \varphi(p^{n-1}) = p^{n-2}(p-1)$. Indeed, equality holds when $M_i$ is cyclic, strict inequality holds when $M_i$ is abelian noncyclic by \cite[Corollary 2.2]{FT23}, and when $M_i$ is nonabelian, Gasch\"utz's theorem \cite{Gas66} gives $|\Aut(M_i)| \ge p|\Inn(M_i)| \ge p^3 \ge p^{n-2}(p-1)$ for $n = 4$, while for $n \ge 5$ we have $|\Aut(M_i)| \ge p^{n-1} > p^{n-2}(p-1)$ (since $|M_i : Z(M_i)| \ge p^2$). Moreover, for $G$ itself ($|G| = p^n$, $n \ge 4$), we have $|\Aut(G)| \ge p^4$: if $|Z(G)| \le p^{n-3}$, then $|\Aut(G)| \ge p|\Inn(G)| \ge p^4$ by \cite{Gas66}; if $|Z(G)| = p^{n-2}$, then $G/Z(G) \cong \Z_p \times \Z_p$ forces $G' \cong \Z_p$ and $\Phi(G) \le Z(G)$, so the group of central automorphisms $\Aut_c(G) \cong \Hom(G/G', Z(G))$ already has order at least $p^{2(n-2)} \ge p^4$.

We distinguish two cases according to the parity of $p$:
\begin{itemize}
    \item If $p \ge 3$, then either $G \cong M_{p^n} = \langle a, b \mid a^{p^{n-1}} = b^p = 1,\, b^{-1}ab = a^{1+p^{n-2}} \rangle$ (which admits a cyclic maximal subgroup and satisfies $|\Aut(M_{p^n})| = p^n(p-1)$, so $f(M_{p^n}) > 1 + (p-1) = p \ge 3 > \frac{5}{2}$), or all $p + 1$ maximal subgroups $M_0, M_1, \dots, M_p$ of $G$ are noncyclic \cite[Theorem 5.4.4]{Gor80}. In the latter case, for $n = 4$ we obtain
    \[
    f(G) > 1 + \frac{(p+1)p^2(p-1) + p^4}{p^4} = 3 - \frac{1}{p^2} \ge \frac{26}{9} > \frac{5}{2}\,.
    \]
    For $n \ge 5$, each noncyclic maximal subgroup $M_i$ ($1 \le i \le p$) has order $p^{n-1} \ge p^4$. If $M_i$ is abelian noncyclic, \cite[Lemma 2.1]{FT23} gives $|\Aut(M_i)| \ge p^{n-2}(p-1)^2 > p^{n-1}$; if $M_i$ is nonabelian, the bound established above yields $|\Aut(M_i)| \ge p^{n-1}$. Thus $|\Aut(M_i)| \ge p^{n-1}$ for all $1 \le i \le p$. Counting the noncyclic subgroups of $M_0$, the $p$ distinct noncyclic maximal subgroups $M_1, \dots, M_p$, and $G$ itself in \eqref{eq:f_noncyclic}, we get
    \[
    f(G) \ge 1 + \frac{\frac{3}{2}p^{n-1} + p \cdot p^{n-1} + p^n}{p^n} = 3 + \frac{3}{2p} > \frac{5}{2}\,.
    \]

    \item If $p = 2$, then either $G$ possesses a cyclic maximal subgroup---in which case $G$ is isomorphic to $D_{2^n}$, $Q_{2^n}$, $SD_{2^n}$, or $M_{2^n}$ \cite[Theorem 5.4.4]{Gor80}, for all of which $f(G) > \frac{5}{2}$ by direct computation---or all maximal subgroups $M_0, M_1, \dots, M_r$ ($r \ge 2$) of $G$ are noncyclic. In the latter case, for $n = 4$ we have $|\Aut(G)| \ge 16$, and hence
    \[
    f(G) \ge \frac{S(M_0)}{16} + \frac{|\Aut(G)|}{16} + \frac{2 \cdot \varphi(8)}{16} \ge \frac{5}{4} + 1 + \frac{1}{2} = \frac{11}{4} > \frac{5}{2}\,.
    \]
    For $n \ge 5$, each noncyclic maximal subgroup $M_i$ ($i = 0, 1, 2$) has order $2^{n-1} \ge 16$, so $|\Aut(M_i)| \ge 2^{n-1}$, while $S(M_0) \ge \frac{5}{2} \cdot 2^{n-1}$ and $|\Aut(G)| \ge 2^{n-1}$. Therefore,
    \[
    f(G) \ge \frac{S(M_0) + |\Aut(M_1)| + |\Aut(M_2)| + |\Aut(G)|}{2^n} \ge \frac{5}{4} + \frac{1}{2} + \frac{1}{2} + \frac{1}{2} = \frac{11}{4} > \frac{5}{2}\,.
    \]
\end{itemize}
This completes the proof.
\end{proof}

We now turn to Problem \ref{prob1}.

\begin{theorem}\label{thm:f_eq_2}
Let $G$ be a finite group. Then $f(G) = 2$ if and only if
\[
G \cong (\Z_3 \rtimes \Z_{2^k}) \times \Z_m,
\]
where $k \in \N^*$, the semidirect product $\Z_3 \rtimes \Z_{2^k} = \langle a, b \mid a^3 = b^{2^k} = 1,\, b^{-1}ab = a^{-1} \rangle$ has nontrivial action, and $m \ge 1$ is an integer with $\gcd(m, 6) = 1$.
\end{theorem}

\begin{proof}
Suppose first that $G \cong (\Z_3 \rtimes \Z_{2^k}) \times \Z_m$ with $\gcd(m, 6) = 1$. Then $\gcd(3 \cdot 2^k, m) = 1$, and by the multiplicativity of $f$ on direct products of coprime orders together with \cite[Proposition 1.3]{FT23}, we obtain
\[
f(G) = f(\Z_3 \rtimes \Z_{2^k}) f(\Z_m) = 2 \cdot 1 = 2.
\]

Conversely, we prove by induction on $|G|$ that for every finite noncyclic group $G$, one has $f(G) \ge 2$, with equality if and only if $G$ is of the stated form.

If $G$ is minimal noncyclic, then \cite[Proposition 1.3]{FT23} gives $f(G) \ge 2$ with equality if and only if $G \cong \Z_3 \rtimes \Z_{2^k}$ for some $k \in \N^*$. Indeed, the only nonabelian minimal noncyclic groups are $Q_8$ (with $f(Q_8) = 4$) and $\Z_p \rtimes \Z_{q^n}$ for distinct primes $p, q$ with $q \mid p-1$, whose sole noncyclic subgroup is $G$ itself, with $|\Aut(G)| = p(p-1)q^{n-1}$, so
\[
f(\Z_p \rtimes \Z_{q^n}) = 1 + \frac{p-1}{q} \ge 2,
\]
with equality if and only if $p = 3$ and $q = 2$. Assume now that $G$ is not minimal noncyclic and that the claim holds for all noncyclic groups of order strictly less than $|G|$. Note that $G$ must be nonabelian, as $f(G) \ge \frac{5}{2} > 2$ for every noncyclic abelian group $G$ by \cite[Theorem 1.4]{FT23}.

Let $M < G$ be a maximal element (with respect to inclusion) among the proper noncyclic subgroups of $G$. By the inductive hypothesis, $f(M) \ge 2$, which means
\[
\sum_{\substack{H \le M \\ H \notin C(M)}} |\Aut(H)| \ge |M|.
\]
Since any subgroup of $G$ containing $M$ is automatically noncyclic, $M$ is in fact a maximal subgroup of $G$. In particular, either $N_G(M) = M$ or $N_G(M) = G$.

\medskip
\noindent\textbf{Case 1: $N_G(M) = M$.} In this case, $M$ has $r = [G : M] \ge 3$ distinct conjugates $M_1, \dots, M_r$ in $G$. Each $M_i \cong M$ contributes at least $|M|$ to the sum in \eqref{eq:f_noncyclic}, and $|\Aut(G)| \ge 4$, so
\[
f(G) \ge 1 + \frac{4 + r|M|}{|G|} \ge 1 + \frac{4 + |G|}{|G|} > 2.
\]

\medskip
\noindent\textbf{Case 2: $N_G(M) = G$, i.e., $M \trianglelefteq G$.} Since $M$ is maximal in $G$, the quotient $G/M$ has no nontrivial proper subgroups, so $[G : M] = p$ is a prime and $|M| = |G|/p$. If $G$ is a $p$-group, Theorem \ref{thm:pgroups} yields $f(G) \ge \frac{5}{2} > 2$. Suppose now that $|G|$ has at least two distinct prime divisors.

\smallskip
\noindent\emph{Subcase 2a: $f(M) > 2$.} Here $\sum_{H \le M,\, H \notin C(M)} |\Aut(H)| > |M|$. If $Z(G) \not\subseteq M$, then $Z(G)$ contains a subgroup $Z$ of order $p$ such that $G = M \times Z$, and Lemma \ref{lem:basic}(c) yields $f(G) \ge f(M) > 2$. Assume henceforth that $Z(G) \le M$. If $Z(G) = 1$, then $|\Aut(G)| \ge |\Inn(G)| = |G| = p|M|$, so counting the noncyclic subgroups inside $M$ together with $G$ itself immediately yields
\[
f(G) > 1 + \frac{|M| + p|M|}{p|M|} = 2 + \frac{1}{p} > 2.
\]
Thus we may also assume $1 < Z(G) \le M$. Choose $x \in G \setminus M$ of $p$-power order so that $G = \langle M, x \rangle$.
\begin{itemize}
    \item Suppose $M$ has a noncyclic Sylow $r$-subgroup $P_r \in \Syl_r(M)$ for some prime $r$ (note that this always occurs when $M$ is nonsolvable, since any finite group with all Sylow subgroups cyclic is solvable).
    
    If $r = p$, then $P = P_p\langle x \rangle \in \Syl_p(G)$ has order $p|P_p| < |G|$, so $f(P) > \frac{5}{2}$ by Theorem \ref{thm:pgroups}. Thus the family $\mathcal{F}$ of noncyclic subgroups of $P$ not contained in $P_p$ contributes at least $\frac{3(p-1)}{2}|P_p| > (p-1)|P_p|$. By the Schur--Zassenhaus theorem \cite{Isa08}, $N_M(P_p) = P_p \rtimes H_{p'}$ for a Hall $p'$-subgroup $H_{p'}$ of $N_M(P_p)$ of order $[N_M(P_p) : P_p]$. Conjugating $\mathcal{F}$ across the $[M : N_M(P_p)]$ Sylow $p$-subgroups of $M$ and extending within $N_G(P_p)$ by the action of $H_{p'}$ (via conjugates $L^y$ or semidirect products $L \rtimes \langle y \rangle$ for $L \in \mathcal{F}$ and $y \in H_{p'}$) multiplies the contribution by $[M : N_M(P_p)] \cdot |H_{p'}| = [M : P_p]$, yielding strictly more than $(p-1)|M|$ outside $M$, so $f(G) > 2$.
    
    If $r = q \ne p$, then $K = P_q\langle x \rangle$ is a proper subgroup of $G$ (if $G = P_q\langle x \rangle$, then $M = P_q \rtimes \Z_{p^{a-1}}$ with $a \ge 2$, and $\Omega_1(P_q)\langle x \rangle < G$ is proper and noncyclic). By the inductive hypothesis, $f(K) \ge 2$, so the noncyclic subgroups of $K$ outside $M$ contribute at least $(p-1)|P_q|$. Applying the same extension via the $[M : N_M(P_q)]$ conjugates of $P_q$ in $M$ and a Hall $q'$-subgroup of $N_M(P_q)$ yields at least $(p-1)|M|$ outside $M$, whence $f(G) > 2$.

    \item Suppose next that all Sylow subgroups of $M$ are cyclic, so $M \cong \Z_u \rtimes \Z_v$ with $\gcd(u,v) = 1$, $u > 1$, $v > 1$. Since $M \trianglelefteq G$, every Sylow subgroup of $\Z_u$ is normal in $G$. Because $f(M) > 2$, either $M$ is not minimal noncyclic---so $\Z_u$ contains a proper nontrivial characteristic subgroup $U_0 \trianglelefteq G$, yielding a proper noncyclic subgroup $U_0 \rtimes (\Z_v\langle x \rangle) < G$ outside $M$ when $x$ acts nontrivially on $U_0$, or $U_0 \times \langle x \rangle$ when $x$ centralizes $U_0$---or $M \cong \Z_q \rtimes \Z_{r^b}$ is minimal noncyclic with $(q,r) \ne (3,2)$ and $r \mid q-1$. In the latter case, $f(M) = 1 + \frac{q-1}{r} \ge 3$, so $\sum_{H \le M,\, H \notin C(M)} |\Aut(H)| \ge 2|M|$. For $p = 2$, since $G \notin C(G)$ and $|\Aut(G)| \ge 4$, we obtain $\sum_{N \notin C(G)} |\Aut(N)| \ge 2|M| + 4 > |G|$ and $f(G) > 2$. For $p \ge 3$, either $p = q$ (so $Q \in \Syl_q(G)$ has order $q^2$ and $G = Q \rtimes \Z_{r^b}$ contains $q$ conjugates of $M$, contributing $q \cdot 2|M| = 2|G| > |G|$) or $p \ne q$ (so $\Z_q\langle x \rangle$ is either $\Z_q \times \Z_p$, giving $p$ conjugates of $M$ or $G = M \times \Z_p$, or $\Z_q \rtimes \Z_p$, whose conjugates and subgroups outside $M$ contribute at least $(p-1)|M|$). In all cases, $f(G) > 2$.
\end{itemize}

\smallskip
\noindent\emph{Subcase 2b: $f(M) = 2$.} By the inductive hypothesis, $M \cong (\Z_3 \rtimes \Z_{2^k}) \times \Z_{m'}$ with $\gcd(m', 6) = 1$. Since $M \trianglelefteq G$ has a unique Hall $\{2,3\}'$-subgroup $\Z_{m'}$, every Sylow subgroup of $\Z_{m'}$ is characteristic in $M$ and hence normal in $G$. Choose $x \in G \setminus M$ of $p$-power order so that $G = \langle M, x \rangle$ and $x^p \in M$.
\begin{itemize}
    \item If $\gcd(p, |M|) = 1$, then $x^p = 1$ and $G = M \rtimes \langle x \rangle \cong M \rtimes \Z_p$. If $\Z_p$ centralizes $M$, then $G \cong (\Z_3 \rtimes \Z_{2^k}) \times \Z_{m'p}$ with $\gcd(m'p, 6) = 1$, giving $f(G) = 2$. If $\Z_p$ acts nontrivially on $M$, then because $|\Aut(\Z_3 \rtimes \Z_{2^k})| = 3 \cdot 2^k$ is coprime to $p$, $\Z_p$ centralizes $\Z_3 \rtimes \Z_{2^k}$ and acts nontrivially on $\Z_{m'}$; hence $G \cong (\Z_3 \rtimes \Z_{2^k}) \times (\Z_{m'} \rtimes \Z_p)$ with $\Z_{m'} \rtimes \Z_p$ noncyclic, so $f(G) = 2 \cdot f(\Z_{m'} \rtimes \Z_p) > 2$.

    \item If $p \mid m'$, then $p \ge 5$. Write $\Z_{m'} = P_0 \times \Z_{m''}$ with $P_0 \in \Syl_p(M)$ cyclic and $\gcd(m'', p) = 1$. Then $P = \langle P_0, x \rangle \in \Syl_p(G)$ is normal in $G$, and by the Schur--Zassenhaus theorem \cite{Isa08}, $G = ((\Z_3 \rtimes \Z_{2^k}) \times \Z_{m''}) \rtimes P$. If $P$ is cyclic and centralizes $(\Z_3 \rtimes \Z_{2^k}) \times \Z_{m''}$, then $G \cong (\Z_3 \rtimes \Z_{2^k}) \times \Z_{m'p}$, giving $f(G) = 2$. If $P$ is noncyclic, then $f(P) \ge \frac{5}{2}$ by Theorem \ref{thm:pgroups} and $f(G) \ge f(\Z_3 \rtimes \Z_{2^k}) f(P) \ge 5 > 2$. If $P$ is cyclic and acts nontrivially on $(\Z_3 \rtimes \Z_{2^k}) \times \Z_{m''}$, it must act nontrivially on $\Z_{m''}$, so $\Z_{m''} \rtimes P$ is noncyclic and $f(G) = 2 \cdot f(\Z_{m''} \rtimes P) > 2$.

    \item If $p = 3$, then $\Z_{m'}$ is a normal Hall $\{2,3\}'$-subgroup of $G$, so $G = H \rtimes \Z_{m'}$ where $|H| = 9 \cdot 2^k$ and $H$ contains $M_0 = \Z_3 \rtimes \Z_{2^k}$ as a normal subgroup of index $3$. Let $P_3 \in \Syl_3(H)$, so $|P_3| = 9$, and let $Q_2 \in \Syl_2(M_0)$ ($Q_2 \cong \Z_{2^k}$), so $H = P_3 \rtimes Q_2$.
    
    If $P_3 \cong \Z_3 \times \Z_3$, then $A = P_3 \cap M_0 \cong \Z_3$ is $Q_2$-invariant. By Maschke's theorem, $P_3 = A \times B$ for a $Q_2$-invariant subgroup $B \cong \Z_3$. If $Q_2$ centralizes $B$, then $H = (A \rtimes Q_2) \times B = M_0 \times B$, and Lemma \ref{lem:basic}(c) together with $P_3 \le H$ ($|\Aut(P_3)| = 48$) gives
    \[
    \sum_{N \notin C(H)} |\Aut(N)| \ge 3 \cdot 2^k + 2(3 \cdot 2^k) + 48 = 9 \cdot 2^k + 48 > |H|.
    \]
    If $Q_2$ acts nontrivially on $B$, then a generator of $Q_2$ inverts both $A$ and $B$, hence inverts every element of $P_3$; consequently, all four subgroups of order $3$ in $P_3$ are $Q_2$-invariant and produce four distinct subgroups of $H$ isomorphic to $M_0$, contributing $4 \cdot (3 \cdot 2^k) = 12 \cdot 2^k > 9 \cdot 2^k = |H|$. Thus $f(G) \ge f(H) > 2$ for all $k \ge 1$.
    
    If $P_3 \cong \Z_9$, then $Q_2 \cong \Z_{2^k}$ acts on $P_3$ by inverting $\Omega_1(P_3) = \Z_3$, hence inverts all of $P_3$. Thus $H \cong \Z_9 \rtimes \Z_{2^k}$ contains $M_0$ (contributing $3 \cdot 2^k$) and $H$ itself ($|\Aut(H)| = 9 \cdot 2^k$), so $\sum_{N \notin C(H)} |\Aut(N)| \ge 12 \cdot 2^k > |H|$, giving $f(G) \ge f(H) > 2$.

    \item If $p = 2$, then $G = H \rtimes \Z_{m'}$ where $|H| = 3 \cdot 2^{k+1}$ and $M_0 = \Z_3 \rtimes \Z_{2^k} \trianglelefteq H$ has index $2$. Let $A \cong \Z_3$ be the unique Sylow $3$-subgroup of $M_0$, so $A \trianglelefteq H$, and let $P_2 \in \Syl_2(H)$ of order $2^{k+1}$, so $H = A \rtimes P_2$. If $P_2 \cong \Z_{2^{k+1}}$ acts nontrivially on $A$, then $H \cong \Z_3 \rtimes \Z_{2^{k+1}}$; if $P_2$ centralizes $\Z_{m'}$, we obtain $G \cong (\Z_3 \rtimes \Z_{2^{k+1}}) \times \Z_{m'}$ with $f(G) = 2$, whereas if $P_2$ acts nontrivially on $\Z_{m'}$, then $\Z_{m'} \rtimes P_2$ is noncyclic and $f(G) > 2$. Note that $P_2 \cong \Z_{2^{k+1}}$ cannot centralize $A$, as $H$ would then be cyclic while $M_0 \le H$ is noncyclic. Finally, if $P_2$ is noncyclic, then $f(P_2) \ge \frac{5}{2}$ by Theorem \ref{thm:pgroups}, while $M_0 \cap P_2 \cong \Z_{2^k}$ is cyclic, so $M_0$ and $P_2$ share no noncyclic subgroups; hence $H$ contains $M_0$ (contributing $3 \cdot 2^k$), $P_2$ (contributing at least $\frac{3}{2} \cdot 2^{k+1} = 3 \cdot 2^k$), and $H$ itself ($|\Aut(H)| \ge 6$), yielding $\sum_{N \notin C(H)} |\Aut(N)| \ge 3 \cdot 2^{k+1} + 6 > |H|$ and $f(G) \ge f(H) > 2$.
\end{itemize}
This completes the proof.
\end{proof}

Finally, we answer Problem \ref{prob3} in the negative by proving that there is no finite group $G$ with $f(G) \in \bigl(2, \frac{5}{2}\bigr)$.

\begin{theorem}\label{thm:gap}
For every finite group $G$, either $f(G) \in \{1, 2\}$ or $f(G) \ge \frac{5}{2}$. In particular,
\[
\mathrm{Im}(f) \cap \left(2, \frac{5}{2}\right) = \emptyset
\]
and therefore $\mathrm{Im}(f)$ is not dense in $[2, \infty)$.
\end{theorem}

\begin{proof}
We proceed by induction on $|G|$. If $G$ is cyclic, then $f(G) = 1$; if $G$ is abelian and noncyclic, then $f(G) \ge \frac{5}{2}$ by \cite[Theorem 1.4]{FT23}. Assume henceforth that $G$ is nonabelian and that the theorem holds for all groups of order strictly smaller than $|G|$.

\medskip
\noindent\textbf{Case 1: $G$ is a minimal noncyclic group.} By \cite[Proposition 1.3]{FT23} and the classification of Miller and Moreno \cite{MM03}, $G$ is isomorphic to $\Z_p \times \Z_p$ ($p$ prime), $Q_8$, or $\Z_p \rtimes \Z_{q^n}$ ($p, q$ distinct primes with $q \mid p-1$). Direct computation yields:
\begin{itemize}
    \item $f(\Z_2 \times \Z_2) = \frac{5}{2}$, and $f(\Z_p \times \Z_p) = 1 + \frac{(p+1)(p-1)^2}{p} \ge \frac{19}{3} > \frac{5}{2}$ for $p \ge 3$;
    \item $f(Q_8) = 4 > \frac{5}{2}$;
    \item $f(\Z_p \rtimes \Z_{q^n}) = 1 + \frac{p-1}{q}$ (since the action has order $q$ and kernel $\Z_{q^{n-1}} = Z(G)$, the only noncyclic subgroup of $G$ is $G$ itself, with $|\Aut(G)| = p(p-1)q^{n-1}$). As $q \mid p-1$, $1 + \frac{p-1}{q}$ is an integer $\ge 2$, so it is either $2$ (for $p = 3, q = 2$) or $\ge 3 > \frac{5}{2}$.
\end{itemize}
Hence the claim holds for all minimal noncyclic groups.

\medskip
\noindent\textbf{Case 2: $G$ is not minimal noncyclic.} Let $M < G$ be a maximal element among the proper noncyclic subgroups of $G$. Since any subgroup containing $M$ is noncyclic, $M$ is a maximal subgroup of $G$, so either $N_G(M) = M$ or $N_G(M) = G$. By the inductive hypothesis, either $f(M) = 2$ or $f(M) \ge \frac{5}{2}$.

\smallskip
\noindent\emph{Subcase 2a: $N_G(M) = M$.} Here $M$ has $r = [G : M] \ge 3$ distinct conjugates $M_1, \dots, M_r$.
\begin{itemize}
    \item If $f(M) \ge \frac{5}{2}$, each $M_i$ contributes $\sum_{H \le M_i,\, H \notin C(M_i)} |\Aut(H)| \ge \frac{3}{2}|M|$ to \eqref{eq:f_noncyclic}, and $|\Aut(G)| \ge 4$, whence
    \[
    f(G) \ge 1 + \frac{4 + r \cdot \frac{3}{2}|M|}{r|M|} = \frac{5}{2} + \frac{4}{|G|} > \frac{5}{2}\,.
    \]

    \item If $f(M) = 2$, then $M \cong (\Z_3 \rtimes \Z_{2^k}) \times \Z_m$ with $k \ge 1$ and $\gcd(m, 6) = 1$ by Theorem \ref{thm:f_eq_2}. For any prime $q \mid m$, the Sylow $q$-subgroup $Q$ of $M$ is characteristic in $M$, so $N_G(M) \le N_G(Q)$; by maximality of $M$, $Q \trianglelefteq G$, and hence $\Z_m \trianglelefteq G$. By the Schur--Zassenhaus theorem \cite{Isa08}, $G = H \rtimes \Z_m$ where $H$ contains $M_0 = \Z_3 \rtimes \Z_{2^k}$ with $N_H(M_0) = M_0$ and $[H : M_0] = r \ge 3$.

    Let $P_2 \in \Syl_2(M_0)$ ($|P_2| = 2^k$). Since $N_H(P_2) \cap M_0 = N_{M_0}(P_2) = P_2$, maximality of $M_0$ in $H$ forces $N_H(P_2) = P_2$ (if $M_0 N_H(P_2) = H$, then $P_2 \trianglelefteq H$, contradicting $M_0 \ne \Z_3 \times P_2$). Thus $P_2 \in \Syl_2(H)$ and $r$ is odd. Let $Q_3 \in \Syl_3(M_0)$ ($|Q_3| = 3$). If $Q_3 \trianglelefteq H$, then $C_H(Q_3) \trianglelefteq H$ has odd order $3r$ and is normalized by $P_2$, so $r$ cannot be composite (otherwise a $P_2$-invariant Sylow subgroup of $C_H(Q_3)$ together with $M_0$ would form an intermediate subgroup); hence $r = q$ is an odd prime.
    \begin{itemize}
        \item If $q = 3$, then $|H| = 9 \cdot 2^k$ and $Q_3 \not\trianglelefteq H$ (otherwise $P_2$ would normalize a Sylow $3$-subgroup $P_3$ of $H$, giving $M_0 < P_3 P_2 < H$). Thus $H$ has $3$ Sylow $3$-subgroups of order $9$, which must be isomorphic to $\Z_3 \times \Z_3$ (if $P_3 \cong \Z_9$, then $Q_3 = \Omega_1(P_3) \trianglelefteq P_3$, so $P_3 \le N_H(Q_3) = M_0$). For $k = 1$, $H \cong (\Z_3 \times \Z_3) \rtimes \Z_2$ has $f(H) = 29 > \frac{5}{2}$; for $k \ge 2$, $H$ contains the $3$ conjugates of $M_0$ ($9 \cdot 2^k$), $P_3$ ($48$), $P_3 \rtimes \Z_{2^{k-1}}$ ($\ge 9 \cdot 2^{k-1}$ by induction), and $H$ itself ($\ge 54$), so $f(G) \ge f(H) > \frac{5}{2}$.

        \item If $q \ge 5$, then $|H| = 3q \cdot 2^k$. If $Q \in \Syl_q(H)$ is normal in $H$, then $H = (\Z_3 \times Q) \rtimes \Z_{2^k}$ contains $q$ conjugates of $M_0$ ($|H|$), the noncyclic subgroup $Q \rtimes \Z_{2^k}$ ($\ge \frac{1}{2}|H|$ by induction), and $H$ itself ($|\Aut(H)| \ge 4$), giving $f(H) > \frac{5}{2}$. Suppose $O_q(H) = 1$. We also have $O_2(H) = 1$: indeed, if $O_2(H) \ne 1$, a minimal normal subgroup $N \le O_2(H)$ of $H$ satisfies $N \le P_2 \cap Z(M_0)$ (as $M_0 N = M_0$), so $k \ge 2$, $|N| = 2$, and $N \le Z(H)$; induction on $k$ applied to $H/N$ (which has $M_0/N \cong \Z_3 \rtimes \Z_{2^{k-1}}$ self-normalizing of index $q$) forces $Q N \trianglelefteq H$, whence $Q \trianglelefteq H$ as $Q \text{ char } Q \times N$, a contradiction. Since $H$ is solvable, its minimal normal subgroup must be $N = O_3(H) \cong \Z_3$. Then $H/N \cong \Z_q \rtimes \Z_{2^k}$ with $q \mid 2^k - 1$ ($k \ge 4$), which again forces $Q \trianglelefteq H$, a contradiction. Thus $f(G) \ge f(H) > \frac{5}{2}$.
    \end{itemize}
\end{itemize}

\smallskip
\noindent\emph{Subcase 2b: $N_G(M) = G$.} Here $[G : M] = p$ is a prime. If $G$ is a $p$-group, Theorem \ref{thm:pgroups} gives $f(G) \ge \frac{5}{2}$. Suppose now that $|G|$ has at least two distinct prime divisors.
\begin{itemize}
    \item If $f(M) \ge \frac{5}{2}$, then $\sum_{H \le M,\, H \notin C(M)} |\Aut(H)| \ge \frac{3}{2}|M|$. If $Z(G) \not\subseteq M$, then $G = M \times Z$ with $|Z| = p$, and Lemma \ref{lem:basic}(c) yields $f(G) \ge f(M) \ge \frac{5}{2}$. Assume now that $Z(G) \le M$. If $Z(G) = 1$, then $|\Aut(G)| \ge |G| = p|M|$, which for $p = 2$ gives
    \[
    f(G) \ge 1 + \frac{\frac{3}{2}|M| + 2|M|}{2|M|} = \frac{11}{4} > \frac{5}{2}\,.
    \]
    In general, choose $x \in G \setminus M$ of $p$-power order so that $G = \langle M, x \rangle$. If $M$ has a noncyclic Sylow $r$-subgroup $P_r \in \Syl_r(M)$ (for $r = p$ or $r = q \ne p$), then $K = P_r\langle x \rangle$ satisfies $f(K) \ge \frac{5}{2}$ (by Theorem \ref{thm:pgroups} when $r = p$, and by the inductive hypothesis when $r = q \ne p$, replacing $K$ by $\Omega_1(P_q)\langle x \rangle$ if $K = G$, since $(q,p) \ne (3,2)$); thus the noncyclic subgroups of $K$ outside $M$ contribute at least $\frac{3(p-1)}{2}|P_r|$, and extending via the $[M : N_M(P_r)]$ conjugates of $P_r$ and a Hall $r'$-subgroup of $N_M(P_r)$ gives at least $\frac{3(p-1)}{2}|M|$ outside $M$, so $f(G) > \frac{5}{2}$. Finally, if all Sylow subgroups of $M$ are cyclic, then $M \cong \Z_u \rtimes \Z_v$ with $\gcd(u,v) = 1$. If some Sylow subgroup of $G$ is noncyclic, it must be the Sylow $p$-subgroup $P$ of $G$, and the same Sylow extension argument gives $f(G) > \frac{5}{2}$. If all Sylow subgroups of $G$ are cyclic, then since $M \trianglelefteq G$ has index $p$ and $Z(G) \le M$, we must have $G \cong \Z_{pu} \rtimes \Z_v$ with $\gcd(pu, v) = 1$ and $\Z_v$ acting nontrivially on the unique subgroup $\Z_p \le \Z_{pu}$; then every noncyclic subgroup $H = \Z_d \rtimes \Z_e \le M$ ($d \mid u$, $e \mid v$) lifts to $\widetilde{H} = \Z_{pd} \rtimes \Z_e \not\subseteq M$ with $|\Aut(\widetilde{H})| \ge \varphi(p)|\Aut(H)| = (p-1)|\Aut(H)|$, contributing at least $\frac{3(p-1)}{2}|M|$ outside $M$ and giving $f(G) \ge \frac{5}{2}$.

    \item If $f(M) = 2$, then $M \cong (\Z_3 \rtimes \Z_{2^k}) \times \Z_m$ with $\gcd(m, 6) = 1$. By Theorem \ref{thm:f_eq_2}, either $f(G) = 2$ or $f(G) > 2$. In the latter case, the classification of the extensions $1 \to M \to G \to \Z_p \to 1$ in the proof of Theorem \ref{thm:f_eq_2} shows that either $G$ contains a noncyclic subgroup $K$ with $f(K) \ge \frac{5}{2}$ (namely $\Z_p \times \Z_p$ when $p \mid m$; $\Z_3 \times \Z_3$ when $p = 3$ and $P_3 \cong \Z_3 \times \Z_3$; $H = \Z_9 \rtimes \Z_{2^k}$ when $p = 3$ and $P_3 \cong \Z_9$, for which $H$ has $3$ conjugates of $\Z_3 \rtimes \Z_{2^k}$ and $|\Aut(H)| = 9 \cdot 2^k$, so $f(H) = 3 > \frac{5}{2}$; or $P_2 \in \Syl_2(G)$ when $p = 2$ and $P_2$ is noncyclic, so $f(P_2) \ge \frac{5}{2}$ by Theorem \ref{thm:pgroups}), or $G \cong (\Z_3 \rtimes \Z_{2^k}) \times (\Z_m \rtimes \Z_p)$ (respectively $(\Z_3 \rtimes \Z_{2^{k+1}}) \times (\Z_m \rtimes \Z_2)$) where $K = \Z_m \rtimes \Z_p$ (respectively $\Z_m \rtimes \Z_2$) is a proper noncyclic subgroup with $f(K) > 2$, hence $f(K) \ge \frac{5}{2}$ by the inductive hypothesis. In all cases, $f(G) \ge f(K) \ge \frac{5}{2}$.
\end{itemize}
Consequently, no finite group $G$ satisfies $2 < f(G) < \frac{5}{2}$\,, which shows that $\mathrm{Im}(f)$ is not dense in $[2, \infty)$.
\end{proof}

\begin{remark}
In fact, $\mathrm{Im}(f)$ is not dense in $\bigl[\frac{5}{2}, \infty\bigr)$ either. For instance, on the class of noncyclic $2$-groups, direct inspection of the proof of Theorem \ref{thm:pgroups} shows that the second smallest value of $f$ after $f(\Z_2 \times \Z_2) = \frac{5}{2}$ is $f(D_8) = \frac{11}{4}$, and more generally $\mathrm{Im}(f) \cap \bigl(\frac{5}{2}, \frac{11}{4}\bigr) = \emptyset$. Table \ref{tab:small_values} illustrates the smallest values of $f(G)$ on noncyclic groups of order at most $16$, highlighting the discrete gaps at the bottom of $\mathrm{Im}(f)$.

\begin{table}[ht]
\centering
\caption{Values of $f(G)$ for representative noncyclic groups of small order.}\label{tab:small_values}
\vspace{4pt}
\begin{tabular}{@{}lccc@{}}
\toprule
Group $G$ & $|G|$ & $S(G)$ & $f(G)$ \\
\midrule
$S_3 \cong \Z_3 \rtimes \Z_2$ & $6$ & $12$ & $2$ \\
$\Z_3 \rtimes \Z_4$ & $12$ & $24$ & $2$ \\
$\Z_2 \times \Z_2$ & $4$ & $10$ & $5/2$ \\
$\Z_4 \times \Z_2$ & $8$ & $22$ & $11/4$ \\
$D_8$ & $8$ & $22$ & $11/4$ \\
$D_{12} \cong S_3 \times \Z_2$ & $12$ & $36$ & $3$ \\
$A_4$ & $12$ & $40$ & $10/3$ \\
$Q_8$ & $8$ & $32$ & $4$ \\
$Q_8 \times \Z_2$ & $16$ & $136$ & $17/2$ \\
\bottomrule
\end{tabular}
\end{table}

It remains an open problem to determine the accumulation points of $\mathrm{Im}(f)$ (or to decide whether $\mathrm{Im}(f)$ is closed and discrete in $[1, \infty)$).
\end{remark}

\bigskip
\noindent{\bf Acknowledgment.} During the preparation of this work, the author used Google Gemini to verify logical steps of the proofs. He reviewed, tested, and validated the output and takes full responsibility for the integrity of the mathematics presented.

\bigskip\noindent{\bf Funding.} The author did not receive support from any organization for the submitted work. 

\bigskip\noindent{\bf Conflict of interest.} The author declares that he has no conflict of interest. 

\bigskip\noindent{\bf Data availability.} The manuscript has no associated data.

\vspace*{3ex}\small

\hfill
\begin{minipage}[t]{5cm}
Marius T\u arn\u auceanu \\
Faculty of Mathematics \\
``Al.I. Cuza'' University \\
Ia\c si, Romania \\
e-mail: {\tt tarnauc@uaic.ro}
\end{minipage}

\end{document}